\pdfoutput=1
\documentclass[journal, 12pt]{IEEEtran}

\usepackage{amsmath,amssymb,amsfonts}

\usepackage{booktabs,multirow}
\usepackage{graphicx}
\usepackage{diagbox}
\usepackage{picins} 
\usepackage{subfigure} 

\usepackage{cite}
\usepackage{url}
\usepackage{color}
\usepackage{multicol}

\newtheorem{assumption}{\bf{Assumption}}
\newtheorem{remark}{\bf{Remark}}
\newtheorem{theorem}{\bf{Theorem}}
\newtheorem{lemma}{\bf{Lemma}}

\newtheorem{definition}{\bf{Definition}}

\begin{document}

\title{Real-time optimal delay minimization algorithms for aircraft on a same runway and dual runways}
\author{Peng Lin, Haopeng Yang
}

\maketitle

\begin{abstract}
In this paper, scheduling problems of aircraft minimizing the total delays on a same runway and on dual runways are studied.
In contrast to the algorithms based on mixed-integer optimization models in existing works, where the optimality and the real-time performance are usually unable to be dealt with at the same time, our work focuses on the interaction mechanism between aircraft coupling with delays and two real-time optimal algorithms are proposed for the four scheduling problems by fully exploiting the combinations of different classes of aircraft based on parallel computing technology.
When $100$ aircraft on dual runways is considered, by using the algorithm in this paper, the optimal solution can be obtained within less than $10$ seconds, while by using the CPLEX software to solve the mix-integer optimization model, the optimal solution cannot be obtained within $1$ hour.
\end{abstract}


\section{Introduction}
With the continuous growth of air transportation, aircraft rescheduling has become more and more important to enhance runway utilization and improve  flight punctuality rate in particular when conflicts arise during aircraft operations or large scale flight delays occur.
Existing studies on aircraft rescheduling can be categorized based on their objective functions: minimizing total delay time\cite{Hu2008,Zhan2010,Balakrishnan2010,Hu2011,Hancerliogullari2013,SAMA2014,Soares2015,LIEDER2015,LIEDER2016,RODRIGUEZDIAZ2017,Brittain2018,Pohl2021,Pohl2022}, minimizing the total deviation of scheduled takeoff and landing times\cite{Beasley2000,Bencheikh2011,SALEHIPOUR2013,FAYE2015,MURCA2015,SABAR2015,Faye2018,SALEHIPOUR2020,Hammouri2020}, minimizing the total delay time for arrival and departure flights\cite{SOLVELING2014,NG2017,SAMA2017,PRAKASH2018,Pang2024}, and minimizing overall operational costs \cite{BENNELL2017,De2018,Gui2023}.
Among these objective functions, minimizing total delay time is the most widely studied due to its significant practical relevance. According to the U.S. Department of Transportation\cite{USDOT2024}, the flight arrival delay rate has increased to $21.90\%$ in 2024, and the number of tarmac delays-defined as domestic flights experiencing delays of more than three hours-rose by $51.2\%$ from 2023 to 2024, leading to substantial economic losses and severe inconvenience for airlines, airports, and passengers.
The effects of flight delays become even more pronounced when airports encounter irregular operations, such as severe weather events or system failures, during which a large number of flights are simultaneously disrupted. In such situations, effective rescheduling and recovery strategies are urgently needed to swiftly and orderly restore airport operations and minimize cascading delays throughout the air transportation network.

Most of the existing works on aircraft rescheduling are based on mixed-integer programming (MIP) and the corresponding resolving algorithms
can generally be categorized into four categories: MIP-based algorithms \cite{Beasley2000,SAMA2014,FAYE2015,Pohl2021,Pohl2022,MURCA2015,PRAKASH2018,Pang2024,Gui2023}, dynamic programming algorithms \cite{Balakrishnan2010,LIEDER2015,LIEDER2016,BENNELL2017,De2018,Faye2018}, heuristic algorithms \cite{SAMA2017,SABAR2015,SALEHIPOUR2020,SOLVELING2014}, and meta-heuristic algorithms \cite{Hu2008, Hu2011,SALEHIPOUR2013, Hancerliogullari2013, RODRIGUEZDIAZ2017, Hammouri2020,Zhan2010,Bencheikh2011,NG2017,Soares2015, Brittain2018}. The biggest problem with these algorithms is that the optimality and the real-time performance cannot be dealt with at the same time. Different from \cite{Hu2008,Zhan2010,Balakrishnan2010,Hu2011,Hancerliogullari2013,SAMA2014,Soares2015,LIEDER2015,LIEDER2016,RODRIGUEZDIAZ2017,Brittain2018,Pohl2021,Pohl2022,Beasley2000,Bencheikh2011,SALEHIPOUR2013,FAYE2015,MURCA2015,SABAR2015,Faye2018,SALEHIPOUR2020,Hammouri2020,SOLVELING2014,NG2017,SAMA2017,PRAKASH2018,Pang2024,BENNELL2017,De2018,Gui2023}, a new theoretical framework was established and real-time optimal algorithms were proposed for
scheduling problem on a same runway and dual runways in \cite{b3b}. But the objective function considered in \cite{b3b} is the minimization of overall operational time.

Following the framework of \cite{b3b}, we consider the total delay time for arrival and departure flights as the objective function, and propose two real-time optimal algorithms for scheduling problem on a same runway and dual runways. Compared to the total operation time, the total delay time need to additionally consider the scheduled operation times to calculate out the delays of all aircraft besides the constraints of time windows and the minimum separation times, which are more complicated to analyze.
Moreover, in contrast to most of the existing works \cite{Hu2008,Zhan2010,Balakrishnan2010,Hu2011,Hancerliogullari2013,SAMA2014,Soares2015,LIEDER2015,LIEDER2016,RODRIGUEZDIAZ2017,Brittain2018,Pohl2021,Pohl2022,Beasley2000,Bencheikh2011,SALEHIPOUR2013,FAYE2015,MURCA2015,SABAR2015,Faye2018,SALEHIPOUR2020,Hammouri2020,SOLVELING2014,NG2017,SAMA2017,PRAKASH2018,Pang2024,BENNELL2017,De2018,Gui2023}, where the obtained solutions were usually not the optimal ones and the errors were also unknown, our algorithms are given based on the analysis of the minimum separation times and the technology of parallel computing technology and can be applied in real-time to find an optimal solution for scheduling problems on a same runway and dual runways. In addition, our algorithms can be used to the RECAT systems, which might have $6$ or more classes of aircraft. This is different from the existing algorithms based on the ICAO separations, where aircraft are usually classified into $3$ classes.

Notations. The operation $A-B$ represents the set that consists of the elements of $A$ which are not elements of $B$; the operation $\phi_0-\phi_1$ represents the sequence which is obtained through modifying the sequence $\phi_0$ by removing the aircraft belonging to the sequence $\phi_1$ and keeping the orders of the rest aircraft unchanged; the symbol $cl_i$ represents the class of aircraft $Tcf_i$; the symbol $/$ represents the meaning of ``or".
\section{Problem formulation}
Without considering other airport constraints, in order to ensure the safety of the aircraft, the minimum separation time between each aircraft and its leading aircraft is only related to their own classes.

Suppose that aircraft can be partitioned into $\eta$ classes in descending order of wake impact, represented as $\mathcal{I}=\{1,2, \cdots, \eta \}$, where $\eta$ is a positive integer, and in general the class $1$ usually refers to $\mathrm{A380}$ aircraft.

Let $P_i$ represent the scheduled takeoff/landing time of aircraft $Tcf_i$. Let $F(\phi, Sr(\phi), Pr_0)=\sum_{i=1}^n h_i(\phi)$ denote the total delays for the aircraft sequence $\phi=\langle Tcf_1, Tcf_2, \cdots, Tcf_n\rangle$ to complete the operation (landing and takeoff) tasks,
where $Sr(\phi)=(S_1(\phi), S_2(\phi), \cdots, S_n(\phi))$ represents the rescheduled takeoff and landing times of aircraft $Tcf_1, Tcf_2, \cdots, Tcf_n$ with $S_1(\phi)\leq S_2(\phi)\leq \cdots \leq S_n(\phi)$, $Pr_0=(P_1, P_2, \cdots, P_n)$, $h_i(\phi)=S_i(\phi)-P_i$ if $S_i(\phi)-P_i>0$ and $h_i(\phi)=0$ if $S_i(\phi)-P_i\leq 0$. For the sake of expression convenience, when no confusion arises, $F(\phi, Sr(\phi), Pr_0)$ can be abbreviated as $F(\phi)$ or $F(\phi, Sr(\phi))$.

The purpose of this paper is to find appropriate operation sequence of aircraft and the corresponding takeoff and landing times to minimize the objective function $F(\phi, Sr(\phi), Pr_0)$, i.e., the total delays of all aircraft, so as to solve the following optimization problem:
\begin{eqnarray}\label{optim1}\begin{array}{lll}\min F(\phi, Sr(\phi), Pr_0)\\
\mbox{Subject to}~~~S_k \in Tf_k=[f_k^{\min}, f_k^{\max}],\\
\hspace {2cm}k =1, \cdots, n,\\
\hspace {2cm}S_{j}-S_i \geq Y_{{i}{j}},~\forall i<j, \\ \hspace{2cm}i ,j=1, \cdots, n,\end{array}\end{eqnarray}
where $Tf_k=[f_k^{\min}, f_k^{\max}]$ represents the set of allowable takeoff or landing times, i.e., the time window constraint, for the aircraft $Tcf_k$ with two constants $f_k^{\min}\leq f_k^{\max}$, $Y_{ij}$ represents the minimum separation time between an aircraft $Tcf_i$ and its trailing aircraft $Tcf_{j}$.

If all of the landing and takeoff operations are interrupted for a period of time, denoted by $[I_{d_1}, I_{d_2}]$ for two constants $I_{d_1}<I_{d_2}$, then the time window constraint $Tf_k$ can be changed into $Tf_k=[f_k^{\min}, f_k^{\max}+I_{d_2}-I_{d_1}]-[I_{d_1}, I_{d_2}]$. For simplicity, we still use $Tf_k=[f_k^{\min}, f_k^{\max}]$ to represent the time window constraint for each aircraft $Tf_k$.

\begin{definition}\label{definition3}{\rm (Relevance) Consider an aircraft sequence $\phi=\langle Tcf_1, Tcf_2, \cdots, Tcf_n\rangle$. Let $S_{ij}=S_j-S_i$ represent the takeoff or landing separation time between an aircraft $Tcf_i$ and its trailing aircraft $Tcf_j$ for $i<j$. If $S_{ij}=Y_{{i}j}$, $i<j$, it is said that aircraft $Tcf_j$ is relevant to the aircraft $Tcf_i$.}\end{definition}

\section{Some necessary assumptions and definitions}

In this section, we give some necessary assumptions and definitions, which are first introduced in \cite{b3b}, to make preparations for the main theoretical results and algorithms.

\subsection{Landing sequences}\label{seclanding}

Let $T_{ij}$ represent the minimum separation time between an aircraft of class $i$ and a trailing aircraft of class $j$ on a single runway without considering the influence of other aircraft.
Let $T_0$ denote the minimum value of all possible separation times between aircraft, which is usually taken as $1$ minute.

Based on the landing separation time standards at Heathrow Airport and the understanding of the physical landing process, we propose the following assumptions.

\begin{assumption}\label{ass1.3.1}{\rm (1) For $i=1,2$, $T_{ii}=1.5T_0$. \\
        (2) For $i=\rho_1, \rho_2$,
        $T_{ii}=T_0+\delta$, where $T_0/8<\delta<T_0/6$ is a positive integer, $\rho_1=3$ and $\rho_2=5$.\\
        (3) For $i\neq 1,2,\rho_1, \rho_2$,
        $T_{ii}=T_0$.}\end{assumption}

\begin{assumption}\label{ass1.3.2}{\rm (1) $T_{21}=1.5T_0$.\\
        (2) For $i>j$ and $i\neq 2$, $T_{ij}=T_0$.}\end{assumption}

\begin{assumption}\label{ass1.3.3}{\rm
(1) For all $i<k<j$, $T_{ik}\leq T_{ij}\leq 3T_0$ and $T_{kj}<T_{ij}\leq 3T_0$.\\
(2) For all $k\leq j\leq i$, $T_{ik}<T_{ij}+T_{jk}$ and $T_{ki}<T_{ji}+T_{kj}$.}\end{assumption}

\begin{definition}\label{definition1}{\rm(Breakpoint aircraft) Consider a landing/takeoff sequence $\phi=\langle Tcf_1, Tcf_2, \cdots, Tcf_{n}\rangle$. If the classes of two consecutive aircraft satisfy that $cl_{i}<cl_{i+1}$, it is said that the aircraft $Tcf_i$ is a breakpoint aircraft of $\phi$.}\end{definition}

\begin{definition}\label{definition2}{\rm(Resident-point aircraft)  Consider a landing or takeoff sequence $\phi=\langle Tcf_1, Tcf_2, \cdots, Tcf_{n}\rangle$. If $S_1>t_0$, $Tcf_1$ is called  a resident-point aircraft of $\phi$, and $S_1-t_0$ is called the resident time of $Tcf_1$. If the aircraft $Tcf_i$ is not relevant to the aircraft $Tcf_{i-1}$, i.e., $S_{(i-1)i}-Y_{(i-1)i}>0$, for $i=2, 3, \cdots, n$,
it is said that $Tcf_i$ is a resident-point of $\phi$ and $S_{(i-1)i}-Y_{(i-1)i}$ is the resident time of $Tcf_i$.
 Consider a mixed landing and takeoff sequence $\phi=\langle Tcf_1, Tcf_2, \cdots, Tcf_{n}\rangle$. Let $\mu_1$ and $\mu_2$ be the largest integers smaller than $i$ such that aircraft $Tcf_{\mu_1}$ is a landing aircraft and aircraft $Tcf_{\mu_2}$ is a takeoff aircraft.
If the aircraft $Tcf_i$ is not relevant to the aircraft $Tcf_{\mu_1}$ and $Tcf_{\mu_2}$, $i=3, 4, \cdots, n$, it is said that $Tcf_i$ is a resident-point aircraft of $\phi$ and $\min\{S_{\mu_1i}(\phi)-Y_{\mu_1{i}}, S_{\mu_2i}(\phi)-Y_{\mu_2{i}}\}$ is the resident time of $Tcf_i$.}\end{definition}

\begin{assumption}\label{ass1.3.4}{\rm (1) For $k=\rho_2$, $T_{(k-1)k}=T_0+\delta$. For $k\neq \rho_2$, $T_{(k-1)k}\geq 1.5T_0$.\\
(2) For all $i\leq k$, when $(i,k)\!\neq\! (\rho_2, \rho_2)$, $T_{(i-1)k}-T_{ik}>2\delta$.\\
(3) For $k=1$, $T_{k(k+1)}>2T_0,$ and for $k=2$, $T_{k(k+1)}>1.5T_0+2\delta$.\\
(4) For $k=1$, all $k+2\leq h \leq \eta$, and all $h\leq j\leq \eta$, $T_{kj}-T_{hj}>0.5T_0$.\\
(5) Let $E=\{\langle3,4\rangle, \langle 3,\eta\rangle\}$ be a sequence set such that
$0.5T_0-\delta\leq T_{2j}-T_{kj}<0.5T_0$ for all $\langle k,j\rangle\in E$ and for all $\langle k,j\rangle\notin E$ with $2<k< j$, $T_{2j}-T_{kj}>0.5T_0$.}\end{assumption}

\begin{definition}{\rm(Class-monotonically-decreasing sequence) If a landing/takeoff sequence $\phi=\langle Tcf_1, Tcf_2, \cdots,$ $ Tcf_n\rangle$ satisfies $cl_1\geq cl_2\geq \cdots\geq cl_n$, then the aircraft sequence $\phi$ is called a class-monotonically-decreasing sequence.
}\end{definition}

\subsection{Takeoff sequences}\label{sectakeoff}

Let $D_{ij}$ represent the minimum separation time between an aircraft of class $i$ and a trailing aircraft of class $j$ on a single runway without considering the influence of other aircraft.

\begin{assumption}\label{ass2.3.1}{\rm (1) When $i=1,2, 3, \eta$, $D_{ii}=(1+1/3)T_0$.
\, (2) When $i\neq 1, 2, 3, \eta$, $D_{ii}=T_0$.}\end{assumption}

\begin{assumption}\label{ass2.3.2}{\rm (1) $D_{21}=(1+1/3)T_0$.
        \, (2) When $i>j$ and $i\geq 3$, $D_{ij}=T_0$.}\end{assumption}

\begin{assumption}\label{ass2.3.3}{\rm
(1) For all $k\leq i\leq j$, $D_{ij}\leq D_{kj}\leq 3T_0$ and $D_{ki}\leq D_{kj}\leq 3T_0$.\\
(2) For all $k\leq j\leq i$, $D_{ik}<D_{ij}+D_{jk}$, $D_{ki}< D_{ji}+D_{kj}$.}\end{assumption}
\begin{assumption}\label{ass2.3.41}{\rm
(1) For $k\!=\!1,\!2$, $D_{k(k+1)}\!\!=\!D_{kk}\!+\tfrac{T_0}{6}$.\\
(2) For $k=3, \rho_2-1$, $D_{k(k+1)}=D_{kk}$.\\
(3) For $k=1$, $D_{k3}=D_{23}+T_0/3$. For $k=1$ and all $k+2\leq j \leq \eta$, $D_{kj}-D_{(k+1)j}=2T_0/3$.\\
(4) For $k=3$, $D_{k\eta}=D_{(k+1)\eta}$, and for $k=\rho_2-1$, $D_{k\eta}=D_{\rho_2\rho_2}+T_0$.\\
(5) For all $2\leq k<j$ and $j\geq 3$ such that $(k,j)\neq (3,\eta)$, $(k,j)\neq (\rho_2-1,\rho_2)$ and $(k,j)\neq (\rho_2-2,\rho_2)$, $D_{kj}=D_{(k+1)j}+T_0/3$.

}\end{assumption}

\subsection{Mixed landing and takeoff sequences on a same runway}\label{seclandingtakeoff}

In this section, the scheduling of mixed takeoffs and landings on a same way for aircraft is discussed. Let $D_T$ denote the minimum  separation time between a landing aircraft and a leading takeoff aircraft. Let $T_D$ denote the minimum separation time between a takeoff aircraft and a leading landing aircraft.

\begin{assumption}\label{ass4.11}{\rm Suppose that $T_0\leq T_D<1.5T_0$ and $T_0\leq D_T<1.5T_0$.}\end{assumption}

\begin{assumption}\label{ass4.21a}{\rm Consider an aircraft sequence $\phi=\langle Tcf_1, Tcf_2, \cdots, Tcf_{n}\rangle$. Suppose that the aircraft $Tcf_{j_1}$ is relevant to $Tcf_{j_0}$, the aircraft $Tcf_{j_2}$ is relevant to $Tcf_{j_1}$, $Y_{j_0j_1}\geq T_D+D_T$ and $Y_{j_1j_2}\geq T_D+D_T$. If the aircraft $Tcf_{j_3}$ is relevant to $Tcf_{j_2}$, then $Y_{j_2j_3}<T_D+D_T$.}\end{assumption}

\begin{definition} {\rm(Path) Consider a sequence $\phi=\langle Tcf_1,$ $Tcf_2, \cdots, Tcf_{n}\rangle$. For any given aircraft $Tcf_i$ and $Tcf_j$, if there exists an aircraft subsequence $\langle \!Tcf_i^0\!, Tcf_i^1\!, \cdots\!, Tcf_i^{\rho}\!\rangle$ for some positive integer $\rho>0$ such that $Tcf_i^0=Tcf_i$, $Tcf_i^\rho=Tcf_j$ and each aircraft $Tcf^h_i$ is relevant to
aircraft $Tcf^{h-1}_i$, $i=1, 2, \cdots, \rho$, then the sequence $\langle Tcf_i^0, Tcf_i^1, \cdots, Tcf_i^{\rho}\rangle$ is said to be a path from the aircraft $Tcf_j$ to the aircraft $Tcf_i$. It is assumed by default that each aircraft has a path to itself.}\end{definition}

\begin{assumption}\label{ass3.4}{\rm Suppose that the aircraft sequence $\phi=\langle Tcf_1, Tcf_2, \cdots, Tcf_{n}\rangle$ is fixed. The aircraft $Tcf_2$ is relevant to aircraft $Tcf_1$, and aircraft $Tcf_{i+2}$ is relevant to aircraft $Tcf_{i+1}$ or aircraft $Tcf_i$, $i=1,2, \cdots, n-2$.}\end{assumption}

\begin{definition}{\rm Consider an aircraft sequence $\langle Tcf_i,\! Tcf_j \!\rangle$. If the aircraft $Tcf_i$ is a takeoff aircraft, and the aircraft $Tcf_j$ is a landing aircraft, it is said that the aircraft sequence forms a takeoff-landing transition at the aircraft $Tcf_i$. If the aircraft $Tcf_i$ is a landing aircraft, and the aircraft $Tcf_j$ is a takeoff aircraft, it is said that the aircraft sequence forms a landing-takeoff transition at the aircraft $Tcf_i$.}\end{definition}
\subsection{Mixed landing and takeoff aircraft on dual runways}\label{secdual}

In this section, scheduling of mixed landing and takeoff on dual runways whose spacing is no larger than $760$ $m$, where all of the landing aircraft lands on one runway and all of the takeoff aircraft take off from the other runway. Let $P_D$ denote the minimum  separation time between a takeoff aircraft and a leading landing aircraft. Let $D_P$ denote the minimum separation time between a landing aircraft and a leading takeoff aircraft.

\begin{assumption}\label{ass4.1}{\rm Suppose that $D_P=T_0$, $P_D=0$, and
for any two aircraft $Tcf_i$ and $Tcf_j$ with the same operation tasks, $Y_{ij}=T_{cl_icl_j}$ or $Y_{ij}=D_{cl_icl_j}$.
}\end{assumption}

Consider a group of landing aircraft and takeoff aircraft operating on dual runways. Let $\Phi=\langle Tcf_1, Tcf_2, \cdots,$ $ Tcf_{n+m}\rangle$ denote the whole mixed landing and takeoff sequence on dual runway, $\Phi_0=\langle Tcf^0_1, Tcf^0_2, \cdots, Tcf^0_n\rangle$ denote the landing aircraft sequence in $\Phi$, $\Phi_1=\langle Tcf^1_{1}, Tcf^1_{2}, \cdots, Tcf^1_{m}\rangle$ denote the takeoff aircraft sequence in $\Phi$. From the definition of aircraft sequence, it follows that $S_1(\Phi)\leq S_2(\Phi)\leq \cdots \leq S_{n+m}(\Phi)$, $S_1(\Phi_0)\leq S_2(\Phi_0)\leq \cdots \leq S_{n}(\Phi_0)$ and $S_1(\Phi_1)\leq S_2(\Phi_1)\leq \cdots \leq S_{m}(\Phi_1)$.

\begin{assumption}\label{ass4.21}{\rm Suppose that for each aircraft $Tcf_i\in \{2,3, \cdots, n+m\}$ in $\Phi$, the aircraft $Tcf_{i}$ is relevant to the aircraft $Tcf_k$ for $k\in \{Tcf_{i-4}, Tcf_{i-3}, Tcf_{i-2}, Tcf_{i-1}\}$. If aircraft $Tcf_{i}$ and $Tcf_{i-1}$ have different (takeoff/landing) operation tasks and aircraft $Tcf_i$ is relevant to $Tcf_{i-1}$ for $i\in \{2,3, \cdots, n+m\}$, there is at least an aircraft $Tcf_{j}\in \{Tcf_1, Tcf_2, \cdots, Tcf_{k-1}\}$ such that $S_{ji}(\Phi)-T_0<Y_{ji}$.}\end{assumption}

\begin{definition}{\rm (T-block, D-block subsequences) Consider a subsequence of the sequence $\Phi$, denoted by $\phi=\langle Tcf_{i_0}, Tcf_{i_0+1}, \cdots, Tcf_{i_1}\rangle$ for two integers $0<i_0<i_1\leq n+m$. Let $\phi_0$ denote the landing aircraft subsequence of $\phi$, and $\phi_1$ denote the takeoff aircraft subsequence of $\phi$. Suppose that there is a path from aircraft $Tcf_{i_1}$ to aircraft $Tcf_{i_0}$,  the operation (landing/takeoff) tasks of $Tcf_{i_0}, Tcf_{i_0+1}$ are different, and $Tcf_{i_0+1}$ is relevant to $Tcf_{i_0}$.

(1) Suppose that $Tcf_{i_1-1}\in \phi_0$, $Tcf_{i_1}\in \phi_1$, for the sequence $\phi_0$, there is a path from its last aircraft $Tcf_{i_1-1}$ to its first aircraft, and for the sequence $\phi_1$, there is a path from its second last aircraft to its first aircraft.
If $Tcf_{i_1}$ is only relevant to $Tcf_{i_1-1}$, it is said that the sequence $\phi$ is a T-block subsequence of $\Phi$, and $(|\phi_1|-1)/(|\phi_0|-1)$ and $S_{i_1}(\phi)-S_{i_2}(\phi)-Y_{i_2i_1}$ are the length and the takeoff time increment of T-block subsequence, where $Tcf_{i_2}$ is the second last takeoff aircraft in $\phi$.

(2) Suppose that $Tcf_{i_1-1}\in \phi_1$, $Tcf_{i_1}\in \phi_0$, for the sequence $\phi_1$, there is a path from its last aircraft $Tcf_{i_1-1}$ to its first aircraft, and for the sequence $\phi_0$, there is a path from the second last aircraft to the first aircraft in $\phi_0$.
If $Tcf_{i_1}$ is only relevant to $Tcf_{i_1-1}$,it is said that the sequence $\phi$ is a D-block subsequence of $\Phi$, and $(|\phi_0|-1)/(|\phi_1|-1)$ and
$S_{i_1}(\phi)-S_{i_2}(\phi)-Y_{i_2i_1}$ are the length and the landing time increment of D-block subsequence, where $Tcf_{i_2}$ is the second last landing aircraft in $\phi$.}\end{definition}

 In fact, under Assumption \ref{ass4.21}, the sequence $\Phi$ can be expressed as a group of block subsequences in the form of, e.g., $\langle TB_1, TB_2, DB_1, TB_3, DB_2, DB_3, TB_4, TB_5, \cdots\rangle$, where each $TB_i$ denotes a T-block subsequence and each $DB_i$ denotes a D-block subsequence. Based on the definitions of block subsequences, we can focus on the switching between T-block subsequences and D-block subsequences to study the time increments of the landing sequence and the takeoff sequence of the whole sequence $\Phi$.

\section{Some theoretical results}

In this section, we present some typical theoretical results which might be frequently used when the algorithms proposed in this paper are implemented.

\begin{theorem}\rm\cite{b3b}\label{theorem1.1} {\rm Consider a landing/takeoff sequence $\phi=\langle Tcf_1, Tcf_2, \cdots, Tcf_{n}\rangle$. Suppose that $Tf_k=[t_0, +\infty]$ for all $k$. Under Assumptions \ref{ass1.3.1}-\ref{ass2.3.41}, if the aircraft $Tcf_j$ is relevant to $Tcf_i$, then $j=i+1$.
}\end{theorem}

\begin{theorem}\rm\cite{b3b}\label{theorem4.11b} {\rm  Consider a mixed landing and takeoff sequence $\phi=\langle Tcf_1, Tcf_2, \cdots, Tcf_{n}\rangle$ on a same runway or dual runways with spacing no larger than 760 m. Suppose that $Tf_j=[t_0, +\infty]$ for all $j$. Under Assumptions \ref{ass1.3.1}-\ref{ass4.21}, the following statements hold.\\
(1) If aircraft $Tcf_j$ is relevant to aircraft $Tcf_i$ and aircraft $Tcf_i$ and $Tcf_j$ have different operation tasks, then $j=i+1$.\\
(2) If aircraft $Tcf_j$ is relevant to aircraft $Tcf_i$ and aircraft $Tcf_i$ and $Tcf_j$ have the same operation tasks, then aircraft $Tcf_{i+1}, Tcf_{i+2}, \cdots, Tcf_{j-1}$ have different operation tasks from aircraft $Tcf_i$ and $Tcf_j$.
}\end{theorem}

\begin{remark}{\rm If an aircraft $Tcf_j$ is relevant to its leading aircraft $Tcf_i$, then the separation time between aircraft $Tcf_i$ and $Tcf_j$ is equal to $Y_{ij}$.}\end{remark}

\begin{theorem}\label{theorem113}{\rm Consider two sequences $\Phi_a=\langle \phi_1, Tcf_{i_0},$ $ Tcf_{i_1}, Tcf_{i_2}, \phi_2, Tcf_{i_3}, Tcf_{i_4}, \phi_3\rangle$, and $\Phi_b=\langle \phi_1, Tcf_{i_0},$ $ Tcf_{i_2}, \phi_2, Tcf_{i_3}, Tcf_{i_1}, Tcf_{i_4}, \phi_3\rangle$, where $\phi_1$, $\phi_2$, $\phi_3$ are three aircraft subsequences containing $m_1,m_2$ and $m_3$ aircraft. Suppose that each aircraft $Tcf_i$ is relevant to its leading aircraft for all $i\neq 1$, $P_i>t_0$ for all $i$, and $S_i(\Phi_a)>P_i$ and $S_i(\Phi_b)>P_i$ for all $i$. Under Assumptions \ref{ass1.3.1}-\ref{ass4.21}, the following statements hold.\\
(1) $F(\Phi_a)-F(\Phi_b)=S_{i_1}(\Phi_a)-S_{i_1}(\Phi_b)+D_1(m_2+2)+D_2(m_3+1)$, where
$D_1=Y_{i_0i_1}+Y_{i_1i_2}-Y_{i_0i_2}$ and $D_2=Y_{i_0i_1}+Y_{i_1i_2}+Y_{i_3i_4}-Y_{i_0i_2}-Y_{i_3i_1}-Y_{i_1i_4}$.\\
(2) If $D_2=0$ and $D_1\leq [S_{i_1}(\Phi_b)-S_{i_1}(\Phi_a)]/(m_2+1)$, $F(\Phi_a)-F(\Phi_b)\leq 0$.
}\end{theorem}

\noindent{Proof:} Theorem \ref{theorem1141}, \ref{theorem11415} and this theorem can be proved by calculations and the detailed proofs are omitted.

 \begin{remark}{\rm Theorem \ref{theorem113}(1) gives a calculation method when the order of an aircraft is adjusted. More detailed discussion can be made according to the specific values of $m_2, m_3, D_1, D_2$ and $S_{i_1}(\Phi_a)-S_{i_1}(\Phi_b)$ which can reduce the computational amount of the algorithms proposed in this paper.}\end{remark}

 \begin{theorem}\label{theorem11419}{\rm Consider a class-monotonically-decreasing landing/takeoff sequence $\Phi_a=\langle Tcf_1, Tcf_2, \cdots, Tcf_n\rangle$. Suppose that each aircraft $Tcf_i$ is relevant to its leading aircraft for all $i\neq 1$,
$S_i(\phi)>P_i>0$ for all $i\in\{1, 2, \cdots, n\}$ and all possible sequences $\phi$. Let $E_1=Y_{12}$ and $E_i=Y_{(i-1)i}+Y_{i(i+1)}-Y_{(i-1)(i+1)}$ for $i\in\{2, \cdots, n-1\}$.  Under Assumptions \ref{ass1.3.1}-\ref{ass4.21}, the following statements hold.\\
(1) Generate a new sequence $\Phi_b$ by moving $Tcf_{j_0}$ to be between $Tcf_{k_0}$ and $Tcf_{k_0+1}$ for $j_0<k_0$. If $E_{j_0}=T_0$, then $F(\Phi_b)\geq F(\Phi_a)$.\\
(2) Suppose that $Y_{k(k+1)}\leq Y_{(k+1)(k+2)}$
for all $k\in\{1,2, \cdots, n-2\}$. Then $(\Phi_a, Sr(\Phi_a))$ is an optimal solution of the optimization problem (\ref{optim1}).
}\end{theorem}
\noindent{Proof:} (1) Since $Y_{ij}\geq T_0$ for all $i,j$ when landing/takeoff aircraft are considered, it follows from Theorem \ref{theorem113} that $F(\Phi_b)\geq F(\Phi_a)$.\\
(2) In the following, we only discuss the case where the number of the aircraft of each class in the sequence $\Phi_a$ is no larger than $1$ and the general case can be discussed in a similar way. Suppose that this statement does not hold. That is, $(\Phi_a, Sr(\Phi_a))$ is not an optimal solution of the optimization problem (\ref{optim1}).
Let $(\Phi_c, Sr(\Phi_c))$ be an optimal solution of the optimization problem (\ref{optim1}), where $\Phi_c=\langle Tcf_{e_1}, Tcf_{e_2}, \cdots, Tcf_{e_n}\rangle$. It follows that $F(\Phi_c)<F(\Phi_a)$. Suppose that $e_1=1, e_2=2, \cdots, e_b=b, e_{b+1}\neq b+1<n$. Note that $\Phi_a$ is a class-monotonically-decreasing landing/takeoff sequence and the number of the aircraft of each class in the sequence $\Phi_a$ is no larger than $1$. Hence, $cl_{e_1}>cl_{e_2}>\cdots>cl_{e_b}>cl_{b+1}>cl_{e_{b+1}}$. It follows that $d>b+1$, $Tcf_{e_{d-1}}$ is a breakpoint aircraft in $\Phi_c$,
and $cl_{e_d}>cl_{e_{d+1}}$, where $e_d=b+1$.
Generate a sequence $\Phi_q=\langle Tcf_{q_1}, Tcf_{q_2}, \cdots, Tcf_{q_n}\rangle$ by moving $Tcf_{b+1}$ to be between $Tcf_{e_b}$ and $Tcf_{e_{b+1}}$ in $\Phi_c$, where $Tcf_{e_k}=Tcf_{q_k}$ for $k=1,2, \cdots, b$.
Note that $Y_{k(k+1)}\leq Y_{(k+1)(k+2)}$.
It is clear that $S_{e_i}(\Phi_c)-S_{q_i}(\Phi_q)=Y_{be_{b+1}}-(Y_{b(b+1)}+Y_{(b+1)e_{b+1}})+Y_{e_{d-1}e_d}+Y_{e_de_{d+1}}-Y_{e_{d-1}e_{d-1}}\geq0$ for $i=d+1,d+2,\cdots,n$ and $S_{e_i}(\Phi_c)-S_{q_i}(\Phi_q)\geq0$ for $i=b+1, b+2, \cdots, d$. From Theorem \ref{theorem113}, it can be obtained that $F(\Phi_d)\leq F(\Phi_c)$. By using a similar approach, it can be proved that $F(\Phi_a)\leq F(\Phi_c)$, which yields a contradiction. Therefore, $(\Phi_a, Sr(\Phi_a))$ is an optimal solution of the optimization problem (\ref{optim1}).

Theorem \ref{theorem11419} gives two rules to judge whether the aircraft orders are optimal. Using similar analysis approaches, we can obtain the following two lemmas.

\begin{lemma}\label{lemma1.31}{\rm  Consider a landing sequence $\Phi_a=\langle\phi_1, \phi_2\rangle$, where $\phi_1=\langle Tcf_1, Tcf_2, \cdots, Tcf_{n_1}\rangle$, $\phi_2=\langle Tcf_{n_1+1},$ $ Tcf_{{n_1}+2}, \cdots, Tcf_{n_2}\rangle$, $\phi_1$ and $\phi_2$ are both class-monotoni- cally-decreasing sequences, $cl_{n_1}<cl_{n_1+1}$, and
$n_1, n_2$ are two positive integers. Merge the aircraft sequences $\phi_1$ and $\phi_2$ to form a new class-monotonically-decreasing sequence $\Phi_b=\langle Tcf_{s_1}, Tcf_{s_2}, \cdots, Tcf_{s_{n_1+n_2}}\rangle$. Suppose that Assumptions \ref{ass1.3.1}-\ref{ass1.3.4} hold for each landing aircraft sequence, each aircraft $Tcf_i$ is relevant to its leading aircraft for all $i\in\{2, 3, \cdots, n_2\}$ and
$S_i(\Phi_a)>P_i>0$ and $S_i(\Phi_b)>P_i>0$ for all $i\in\{1, 2, \cdots, n_2\}$, the following statements hold.\\
(1) If $(cl_i, cl_{n_1}, cl_{n_1+1})=(\rho_2, \rho_2-1, \rho_2)$ for some $i\in \{1, 2, \cdots, n_1-1\}$, $F(\Phi_a)<F(\Phi_b)$.\\
(2) If $cl_{n_1}\notin\{1,2\}$ and $(cl_i, cl_{n_1}, cl_{n_1+1})\neq(\rho_2, \rho_2-1, \rho_2)$ for all $i\in \{1, 2, \cdots, n_1-1\}$, $F(\Phi_a)-F(\Phi_b)>0$.\\
(3) If $cl_{n_1}\in\{1,2\}$, $F(\Phi_a)-F(\Phi_b)>0$.}\end{lemma}

\begin{lemma}\label{lemma2.31}{\rm  Consider a  takeoff sequence $\Phi_a=\langle\phi_1, \phi_2\rangle$, where $\phi_1=\langle Tcf_1, Tcf_2, \cdots, Tcf_{n_1}\rangle$, $\phi_2=\langle Tcf_{n_1+1},$ $ Tcf_{{n_1}+2}, \cdots, Tcf_{n_2}\rangle$, $\phi_1$ and $\phi_2$ are both class-monotoni- cally-decreasing sequences, $cl_{n_1}<cl_{n_1+1}$, and
$n_1>1$ and $n_2>1$ are two positive integers. Merge the aircraft sequences $\phi_1$ and $\phi_2$ to form a new class-monotonically-decreasing sequence $\Phi_b$. Suppose that Assumptions \ref{ass2.3.1}-\ref{ass2.3.41} hold for each takeoff aircraft sequence, each aircraft $Tcf_i$ is relevant to its leading aircraft for all $i\in\{2, 3, \cdots, n_2\}$ and
$S_i(\Phi_a)>P_i>0$ and $S_i(\Phi_b)>P_i>0$ for all $i\in\{1, 2, \cdots, n_2\}$, the following statements hold.\\
(1) Suppose that $(cl_{n_1}, cl_{n_1+1})=(\rho_2-1, \rho_2)$. Then $F(\Phi_a)-F(\Phi_b)=0$.\\
(2) Suppose that $(cl_{n_1}, cl_{n_1+1}, j)=(3,4,3)$ for some $j\in\{n_1+2,n_1+3,\cdots, n_1+n_2\}$. Then $F(\Phi_a)>F(\Phi_b)$.\\
(3) Suppose that $(cl_{n_1}, cl_{n_1+1}, j)=(\eta-1,\eta,\eta)$ for some $j\in\{1, 2, \cdots, n_1-1\}$. Then $F(\Phi_a)<F(\Phi_b)$.\\
(4) Suppose that $cl_{n_2}\leq 2$ and $(cl_{n_1}, cl_{n_1+1},j)=(2,3,3)$ for some $j\in\{1, 2, \cdots, n_1-1\}$. Then {$F(\Phi_a)=F(\Phi_b)$}.

}\end{lemma}

\begin{theorem}\label{theorem1141}{\rm Consider a landing/takeoff sequence $\Phi_a=\langle Tcf_1, Tcf_2, \cdots, Tcf_n\rangle$ on a same runway. Generate a new sequence $\Phi_b$ by moving $Tcf_{j_0}$ to be between $Tcf_{k_0}$ and $Tcf_{k_0+1}$ for $j_0>k_0+1$.  Suppose that $S_1(\Phi_a)=S_1(\Phi_b)=t_0$ and each aircraft $Tcf_i$ is relevant to its leading aircraft for all $i\neq 1$ in $\Phi_a$ and $\Phi_b$. Under Assumptions \ref{ass1.3.1}-\ref{ass4.21}, the following statements hold.\\
(1) Suppose that $S_{j_0}(\Phi_a)\leq P_{j_0}$. It follows that $F(\Phi_b)\geq F(\Phi_a)$. If $P_{j}<S_{j}(\Phi_b)$ for some $j\in \{k_0+1, k_0+2, \cdots, j_0-1\}$, then $F(\Phi_b)> F(\Phi_a)$.\\
(2) Suppose that $S_{j_0}(\Phi_b)>P_{j_0}$, $Y_{j_0(k_0+1)}+Y_{k_0j_0}-Y_{k_0(k_0+1)}=T_0$ and $Y_{j(j+1)}\geq T_0$ for all $j\in\{1, 2, \cdots, n-1\}$. It follows that $F(\Phi_b)\leq F(\Phi_a)$. Further, if $Y_{j_0(j_0+1)}+Y_{(j_0-1)j_0}-Y_{(j_0-1)(j_0+1)}>T_0$ or $Y_{j(j+1)}>T_0$ for some $j\in \{k_0+1, k_0+2, \cdots, j_0-2\}$, then $F(\Phi_b)<F(\Phi_a)$.
}\end{theorem}

\begin{remark}{\rm Theorem \ref{theorem1141}(1) shows that if the forward movement of an aircraft without delay in sequences would result in nondecreasement of the objective function $F(\cdot)$.}\end{remark}

 \begin{theorem}\label{theorem14}{\rm Consider a sequence $\Phi_a=\langle Tcf_1, Tcf_2, \cdots,$ $Tcf_n\rangle$ on a same runway or dual runways with spacing no larger than 760 m. Suppose that $S_k(\phi)>P_k>0$ for all $k\in\{1, 2, \cdots, i\}$ and all possible sequences $\phi$.
and $(\phi_{i}, Sr(\phi_{i}))$ is one solution of the optimization problem (\ref{optim1}), where $\phi_{i}=\langle Tcf_1, Tcf_2, \cdots, Tcf_i\rangle$. Suppose that each aircraft $Tcf_i$ is relevant to its leading aircraft for all $i\neq 1$ and $i\neq j_0+1$.
If ${S}_j(\phi_{i})+T_{\sigma}\in [f_{j}^{\min}, f_{j}^{\max}]$ for all $j\in\{j_0,j_0+1,\cdots,i\}$ and the operation times of $Tcf_{j_0-1}$ and $Tcf_{j_0}$ are fixed as $S_{j_0-1}(\phi_{i})$ and $S_{j_0}(\phi_{i})+T_{\sigma}$ for a constant $T_{\sigma}>0$, then the function of $F(\phi, {Sr}(\phi),Pr_0)-S_{(j_0-1)j_0}(\phi)$ can be minimized by $(\phi_{i}, \overline{Sr}(\phi_{i}))$, where $\overline{Sr}(\phi_{i})=[S_1(\phi_{i}), \cdots, S_{j_0-2}(\phi_{i}), \bar{S}_{j_0-1}(\phi_{i}), \bar{S}_{j_0}(\phi_{i}), {S}_{j_0+1}(\phi_{i})+T_{\sigma}, \cdots, {S}_{i}(\phi_{i})+T_{\sigma}]$.
}\end{theorem}

\noindent{Proof:} For landing/takeoff sequences, each aircraft is relevant to only its leading aircraft. Since $\bar{S}_j(\phi_{i})\in [f_{j}^{\min}, f_{j}^{\max}]$ for all $j\in\{j_0,j_0+1,\cdots,i\}$ and $(\phi_{i}, Sr(\phi_{i}))$ is one solution of the optimization problem (\ref{optim1}) for aircrafts $Tcf_{11},Tcf_{12}, \cdots, Tcf_{1i}$, the function of $F(\phi, \overline{Sr}(\phi))-S_{(j_0-1)j_0}(\phi)$ can be minimized by $(\phi_{i}, \overline{Sr}(\phi_{i}))$.

This theorem actually studies a preservation problem of the sequence optimality when some time interval is occupied by some operations, e.g.,
new aircraft are inserted. Specifically, Theorem \ref{theorem14} discusses the case when the optimality can be maintained.
This theorem is very useful to analyze the optimality of the sequences to reduce the computation of algorithms when new aircrafts are inserted into the original sequences while keeping the orders of its adjacent aircraft unchanged.

 \begin{theorem}\label{theorem11415}{\rm Consider a mixed landing and takeoff sequence $\Phi_a=\langle \phi_{a1}, \phi_{a2}\rangle$ on a same runway or dual runways, where $\phi_{a1}=\langle Tcf_1, Tcf_2, \cdots, Tcf_{m_0}\rangle$ for some integer $m_0>0$, $\phi_{a2}=\langle Tcf_{m_0+1}, Tcf_{m_0+2}, \cdots, Tcf_{m_0+m_1}\rangle$ for some integer $m_1>0$ and $\phi_{a2}$ contains $m_2$ delayed aircraft.
Let $\Phi_b=\langle \phi_{b1}, \phi_{a2}\rangle$ be a sequence with the same group of aircraft as $\Phi_a$, where $\phi_{a2}$ contains at least $m_2$ delayed aircraft. Suppose that
each aircraft $Tcf_i$ is relevant to its leading aircraft for all $i\neq 1$ in $\Phi_a$ and $\Phi_b$. Under Assumptions \ref{ass1.3.1}-\ref{ass4.21},
if $m_2[S_{m_0+1}(\Phi_{b})-S_{m_0+1}(\Phi_{a})]+F(\phi_{b1}, Sr(\phi_{b1}))-F(\phi_{a1}, Sr(\phi_{a1}))>0$, then $F(\Phi_a)-F(\Phi_b)<0$.
}\end{theorem}

Under this theorem, when the number of delayed aircraft in $\phi_{a2}$ is large, by simple calculations similar to Theorem \ref{theorem113}, it can be obtained that the operation time of the first aircraft of $\phi_{a2}$ should be minimized when $F(\phi_{b1}, Sr(\phi_{b1}))-F(\phi_{a1}, Sr(\phi_{a1}))$ is small, which might be equivalent to be a makespan problem and might be solved by using the algorithm in \cite{b3b}.

\begin{theorem} \rm\cite{b3b} {\rm Consider a mixed landing and takeoff sequence $\Phi_a=\langle Tcf_1, Tcf_2, Tcf_3, Tcf_4\rangle$ on a same runway, where $S_{12}(\Phi_a)\geq Y_{12}$, aircraft $Tcf_3$ is relevant to aircraft $Tcf_1$ or $Tcf_2$ and aircraft $Tcf_4$ is relevant to aircraft $Tcf_2$ or $Tcf_3$. Generate a new sequence $\Phi_b$ by inserting an aircraft $Tcf_5$ to be between aircraft $Tcf_2$ and $Tcf_3$, where $S_{1}(\Phi_a)=S_{1}(\Phi_b)$, $S_{2}(\Phi_a)=S_{2}(\Phi_b)$, aircraft $Tcf_5$ is relevant to aircraft $Tcf_1$ or $Tcf_2$, aircraft $Tcf_3$ is relevant to aircraft $Tcf_2$ or $Tcf_5$ and aircraft $Tcf_4$ is relevant to aircraft $Tcf_5$ or $Tcf_3$. Let $\omega_{min}(\Phi_a, Tcf_5)=\min\{S_3(\Phi_b)-S_3(\Phi_a),S_4(\Phi_b)-S_4(\Phi_a)\}$. The following statements hold.\\
(1) Suppose that aircraft $Tcf_2$, $Tcf_3$ and $Tcf_5$ are all landing (takeoff) aircraft. It follows that $\omega_{min}(\Phi_a, Tcf_5)=Y_{25}+Y_{53}-Y_{23}$.\\
(2) Suppose that aircraft $Tcf_2$ and $Tcf_5$ are both landing (takeoff) aircraft and $Tcf_3$ is a takeoff (landing) aircraft.
It follows that $\omega_{min}(\Phi_a, Tcf_5)\geq\min\{T_D+D_T-Y_{13},T_D+D_T-Y_{24},0\}+Y_{25}$.\\
(3) Suppose that aircraft $Tcf_3$ and $Tcf_5$ are both landing (takeoff) aircraft and $Tcf_2$ is a takeoff (landing) aircraft. It follows that $\omega_{min}(\Phi_a, Tcf_5)\geq\min\{T_D+D_T-Y_{13},T_D+D_T-Y_{24},0\}+Y_{53}$.\\
(4)  Suppose that aircraft $Tcf_2$ and $Tcf_3$ are both landing (takeoff) aircraft and $Tcf_5$ is a takeoff (landing) aircraft. It follows that $\omega_{min}(\Phi_a, Tcf_5)\geq \max\{D_T+T_D-Y_{23},0\}$.
}\end{theorem}

\begin{theorem}{\rm Consider a mixed landing and takeoff sequence $\Phi=\langle Tcf_1, Tcf_2, \cdots, Tcf_{n+m}\rangle$ on dual runways with spacing no larger than $760$ $m$. Suppose that $Tcf_i, Tcf_{i+1}, Tcf_{i+2}, \cdots, Tcf_{n+m}$ are the last $n+m-i+1$ aircraft of $\Phi$, aircraft $Tcf_{i+1}$ and $Tcf_{i+2}$ are landing aircraft, $Tcf_i, Tcf_{i+3}, Tcf_{i+4}, \cdots, Tcf_{n+m}$ are takeoff aircrafts, and $Tcf_i, Tcf_{i+1}, Tcf_{i+2}, Tcf_{i+3}$ are part of a T-block subsequence where $Tcf_{i+3}$ is not relevant to $Tcf_i$. Generate a new sequence $\Phi_a$ by exchanging the orders of the aircraft $Tcf_{i+2}$ and $Tcf_{i+3}$. Suppose that aircraft $Tcf_{i+3}$ is relevant to aircraft $Tcf_{i+2}$ in $\Phi$, aircraft $Tcf_{i+3}$ is relevant to aircraft $Tcf_{i+1}$ in $\Phi_a$ and each $Tcf_k$ is relevant to aircraft $Tcf_{k-1}$ for all $k=i+4, \cdots, n+m$ in $\Phi$ and $\Phi_a$.
If $(n+m-i-2)(S_{i+3}(\Phi)-S_{i+3}(\Phi_a))+S_{i+2}(\Phi)-S_{i+2}(\Phi_a)>0$, then $F(\Phi)-F(\Phi_a)>0$.}\end{theorem}

This theorem is important and can be used to determine whether the sequence is optimal for the optimization problem (\ref{optim1}) by analyzing the aircraft at the tail of the sequence.

\section{Algorithms}
In this section, we propose algorithms to find the optimal solution for the optimization problem (\ref{optim1}) based on the theoretical results given in \cite{b3b} and this paper.

\subsection{Algorithm 1}
\textbf{Algorithm 1.} Suppose that there are totally $n$ aircraft, denoted by $Tcf_{01}, Tcf_{02}, \cdots, Tcf_{0n}$, and there is at least a feasible sequence for them to land/takeoff without conflicts.

Step 1. Arrange the aircraft in ascending order of their scheduled landing/takeoff times, denoted by $\langle Tcf_{11}, Tcf_{12}, \cdots, Tcf_{1n}\rangle$.

Step 2. Search for the optimal sequence $\phi_2$ for the optimization problem (\ref{optim1}).

Step 3. Suppose that $\phi_2=\langle Tcf_1, Tcf_2\rangle$ is an optimal sequence of the optimization problem (\ref{optim1}) for aircraft $Tcf_1$ and $Tcf_2$. Search for the optimal sequence of the optimization problem (\ref{optim1}) for aircraft $Tcf_1$, $Tcf_2$ and $Tcf_{13}$.

Step $i$, $i=4,5, \cdots,n$. Suppose that $\phi_{i-1}=\langle Tcf_1,Tcf_2,$ $ \cdots, Tcf_{i-1}\rangle$ is an optimal sequence of the optimization problem (\ref{optim1}) for aircraft $Tcf_1, Tcf_2, \cdots, Tcf_{i-1}$.
Search for an optimal sequence of the optimization problem (\ref{optim1}) for aircraft $Tcf_1, Tcf_2, \cdots, Tcf_{i}$. as follows.

$(i-1)$. Insert aircraft $Tcf_{1i}$ between any two adjacent aircraft in $\phi_{i-1}$ under the time window constraints without changing the orders of the aircraft in $\phi_{i-1}$. Let $f^i_{inc}$ be the smallest value of the objective functions $F(\cdot)$ among all the generated new sequences and $\phi^i_{inc}$ denote the corresponding sequence.

$(i-2)$.  Construct a sequence set $F_{inc}=\{\langle h_1, h_2, 1i, h_3, $ $h_4\rangle \mid h_1, h_2, h_3, h_4\in \{1,2,\cdots,i-1\}\}$. Note here that since the minimum separation times are related to only the classes and the operation tasks of the aircraft, we can classify the aircraft to form the set $F_{inc}$ according to the classes and the operation tasks of the aircraft.

$(i-3)$. According to the elements of the set $F_{inc}$, adjust the aircraft orders in the sequence $\phi_{i-1}$ and insert the aircraft $Tcf_{1i}$ to generate a new sequence $\bar{\phi}_{i}$ under the condition $F(\bar{\phi}_{i})<f^i_{inc}$ such that $\bar{\phi}_{i}$ contains the subsequence $\langle Tcf_{h_1}, Tcf_{h_2}, Tcf_{1i}, Tcf_{h_3}, Tcf_{h_4}\rangle$ and $\langle h_1, h_2, h_3, h_4\rangle\in F_{inc}$. Minimize $F(\bar{\phi}_{i})$ based on the obtained theoretical results, and compare the values of all possible $F(\bar{\phi}_{i})$ so as to find
the optimal sequence of aircraft $Tcf_{11},Tcf_{12}, \cdots, Tcf_{1i}$.

When landing/takeoff sequences are considered and all aircraft are delayed in all possible sequences, if the aircraft $Tcf_{1i}$ is the $k$th aircraft of $\bar{\phi}_{i}$ and
the inequality is not satisfied\begin{eqnarray}\label{eaf1}\begin{array}{lll}\omega_{min}(\langle Tcf_{h_1}, Tcf_{h_2}, Tcf_{h_3}, Tcf_{h_4}\rangle, Tcf_{1i})\\\times(i-k)+h_{1i}(\bar{\phi}_{i})< f^i_{inc}-F(\phi_{i-1})\end{array}\end{eqnarray} for any $\langle h_1, h_2, h_3, h_4\rangle\in F_{inc}$, then $F(\bar{\phi}_{i})>f^i_{inc}$. This fact can be used to reduce the computation amount of the algorithm.

In step $i.(i-3)$ of Algorithm 1, we can
 discuss the sequence optimality in the following cases and search for the optimal sequence by comparing the values of the objective function $F(\cdot)$ of all possible sequences based on parallel computing technology.

Case 1.1 The number of breakpoint aircraft in $\phi_{i}$ is no larger than that in $\phi_{i-1}$ for landing/takeoff sequences.\\
Case 1.2.  The number of breakpoint aircraft in $\phi_{i}$ is larger than that in $\phi_{i-1}$ for landing/takeoff sequences.\\
Case 2.1.  The number of landing-takeoff and takeoff-landing transitions in $\phi_{i}$ is no larger than that in $\phi_{i-1}$ for mixed landing and takeoff sequences on a same runway.\\
Subcase 2.1.1. The number of breakpoint aircraft in $\phi_{i}$ is no larger than that in $\phi_{i-1}$.\\
Subcase 2.1.2. The number of breakpoint aircraft in $\phi_{i}$ is larger than that in $\phi_{i-1}$.\\
Case 2.2.  The number of breakpoint aircraft in $\phi^0_{i}$ is larger than that in $\phi_{i-1}$ for mixed landing and takeoff sequences on a same runway.\\
Subcase 2.2.1. The number of breakpoint aircraft in $\phi_{i}$ is no larger than that in $\phi_{i-1}$.\\
Subcase 2.2.2. The number of breakpoint aircraft in $\phi_{i}$ is larger than that in $\phi_{i-1}$.

When the number of breakpoint aircraft in $\phi_{i}$ is 2 larger than that in $\phi_{i-1}$, $\phi_{i}$ is usually not the optimal sequence when the aircraft subsequence $\langle Tcf_i, Tcf_j\rangle$ with $cl_i=\rho_2-1$ and $cl_j=\rho_2$ is not involved.

Note that $\phi_{i-1}$ might be composed of multiple class-monotonically-decreasing subsequences $\phi^1, \phi^2, \cdots, \phi^{s}$. We can consider the insertion of $Tcf_{1i}$ on each subsequence $\phi^k$ with some corresponding aircraft order adjustments so as to reduce the complexity of the analysis of the optimality of the sequence.

Note that each aircraft might be relevant to its leading aircraft and nearest aircraft ahead with the same operation task for mixed landing and takeoff sequences on a same runway. The insertion of the aircraft $Tcf_{1i}$ might have direct impact on the following two aircraft and be affected by the preceding two aircraft. So, when the aircraft $Tcf_{1i}$ is inserted, we might need to consider all possible combinations of $5$ consecutive aircraft of the form $\langle Tcf_{h_1}, Tcf_{h_2}, Tcf_{1i}, Tcf_{h_3}, Tcf_{h_4}\rangle$ in the new generated sequence.

From Theorem \ref{theorem1141}, when the number of delayed aircraft in the considered sequence is large, the operation time of the first half of the sequence should be minimized for the optimization problem (\ref{optim1}). Based on the obtained theoretical results, especially Theorem \ref{theorem1141}, we can simplify the generation of the combinations of $5$ consecutive aircraft of the form $\langle Tcf_{h_1}, Tcf_{h_2}, Tcf_{1i}, Tcf_{h_3}, Tcf_{h_4}\rangle$ to reduce the computation complexity. Moreover, since aircraft might be relevant to only its leading aircraft, then the combination of $\langle Tcf_{h_1}, Tcf_{h_2}, Tcf_{1i}, Tcf_{h_3}, Tcf_{h_4}\rangle$ can be simplified as $\langle Tcf_{h_2}, Tcf_{1i}, Tcf_{h_3}\rangle$ in particular for landing/takeoff sequences.

Due to the constraints of time windows of aircraft, breakpoint aircraft and resident-point aircraft might be generated which might increase the value of the objective function $F(\cdot)$. To analyze the optimal sequence of the aircraft, the emphasis should be imposed on the breakpoint aircraft and the resident-point aircraft as well as the separation times between the consecutive aircraft of the same class that are larger than $T_0$ by trying not to increase their numbers. Moreover, when the orders of aircraft need to be adjusted, it is better to fully consider the classes of the aircraft rather than the aircraft themselves.

\subsection{Algorithm 2}

In the following, we propose an algorithm to deal with the scheduling problem of landing and takeoff aircraft on dual runways with spacing no larger than $760$ $m$.

\textbf{Algorithm 2.} Suppose that there are totally $n+m$ aircraft composed of $n$ landing aircraft and $m$ takeoff aircraft, denoted by $Tcf_{1}, Tcf_{2}, \cdots, Tcf_{n+m}$, and there is at least a feasible sequence for them to land/takeoff without conflicts. Let $Tcf^0_1, Tcf^0_2, \cdots, Tcf^0_n$ denote all the landing aircraft, and $Tcf^1_{1}, Tcf^1_{2}, \cdots, Tcf^1_{m}$ denote all the takeoff aircraft.
The main steps are as follows.

Step 1. Use Algorithm 1 to find the optimal sequence, denoted by $\Phi_a$, for the landing aircraft to minimize $F(\Phi_a)$,
 and the optimal sequence for the takeoff aircraft, denoted by $\Phi_b$ to minimize $F(\Phi_b)$. Suppose that $(\Phi, Sr(\Phi))$ is an optimal solution of the optimization problem (\ref{optim1}). It is clear that $F(\Phi)\geq F(\Phi_a)+F(\Phi_b)$.

Step 2. (2.1) According to the landing/takeoff time increments of block subsequences, classify all possible block subsequences into several sets, which can be processed offline.

It should be noted that all possible block subsequences can be obtained by an enumeration method. It should also be noted that when the number of breakpoint aircraft in a block subsequence is fixed, the total number of all possible block subsequences is not large under the minimum separation standards of Heathrow Airport and the RECAT-EU system. This observation can be used to reduce the computational burden in Step 5 of the algorithm.

(2.2) Focus on the combinations of different block subsequences and analyze the landing/takeoff time increments between any two consecutive blocks, where the last two aircraft of the first block are identical to the first two aircraft of the second block. According to the landing/takeoff time increments of two consecutive block subsequences, classify all possible combinations of two consecutive subsequences into several sets.

Step 3. (3.1) Match the takeoff aircraft with the landing sequence $\Phi_a$ to generate a sequence, denoted by $\Phi_c=\langle Tcf_{c1}, Tcf_{c2}, \cdots, Tcf_{c(n+m)}\rangle$, minimizing $F(\Phi_c)$ under the condition that $F(\Phi_{LA}^c, Sr_{LA}(\Phi_c))=F(\Phi_a)$, where $\Phi^c_{LA}$ and $Sr_{LA}(\Phi_c)$ denote the landing sequence and the corresponding landing time vector in $\Phi_c$.

(3.2) Let $c_z$ denote the number of the aircraft of $\Phi_c$ which take off in $(S_n(\Phi_a),+\infty)$. Insert the last $c_z$ aircraft of $\Phi_c$ into the subsequence of $\Phi_c$, $\Phi^{sub}_{c}=\langle Tcf_{c1}, Tcf_{c2}, \cdots, Tcf_{c(n+m-c_z)}\rangle$, to generate a new sequence $\Phi_{c3}$ and minimize $F(\Phi_{c3})$ without changing the orders of the landing aircraft based on the sets defined Step 2. Note here that the sequence $\Phi_{c3}$ can contain D-block subsequences.

It should be noted that the optimal value of the objective function $F(\cdot)$ for the optimization problem (\ref{optim1}) lies in $[F(\Phi_a)+F(\Phi_b), F(\Phi_{c3})]$, and the final optimal sequence usually contains at least $c_z$ D-block subsequences, which can be used as an initial condition for the search in step 4.

Step 4. (4.1) Match the landing aircraft with the takeoff sequence $\Phi_b$ to generate a new sequence, denoted by $\Phi_d=\langle Tcf_{d1}, Tcf_{d2}, \cdots, Tcf_{d(n+m)}\rangle$, and optimize $\Phi_d$ such that $F(\Phi^d_{TA}, Sr_{TA}(\Phi_d))=F(\Phi_b)$, where $\Phi^d_{TA}$ and $Sr_{TA}(\Phi_d)$ denote the takeoff sequence and the corresponding takeoff time vector in $\Phi_d$.

(4.2) Let $d_z$ denote the number of the aircraft of $\Phi_d$ which land in $(S_m(\Phi_b),+\infty)$. Insert the last $d_z$ aircraft of $\Phi_d$ into the subsequence of $\Phi_d$, $\Phi^{sub}_{d}=\langle Tcf_{d1}, Tcf_{d2}, \cdots, Tcf_{d(n+m-d_z)}\rangle$, to generate a new sequence $\Phi_{d3}$ and minimize $F(\Phi_{d3})$ without changing the orders of the takeoff aircraft based on the sets defined Step 2. Note here that the sequence $\Phi_{d3}$ can contain T-block subsequences.

It should be noted that the optimal value of the objective function $F(\cdot)$ for the optimization problem (\ref{optim1}) lies in $[F(\Phi_a)+F(\Phi_b), \min\{F(\Phi_{c3}), F(\Phi_{d3})\}]$.

Step 5. The main idea of this step is to implement a full space search for an optimal solution $(\Phi_{c4}, Sr(\Phi_{c4}))$ for the optimization problem (\ref{optim1}) under the condition that $F(\Phi_a)+F(\Phi_b)\leq F(\Phi_{c4})\leq \min\{F(\Phi_{c3}), F(\Phi_{d3})\}$ starting from the initial condition that $\Phi^{c4}_{LA}=\Phi_{a}$ or $\Phi^{c4}_{TA}=\Phi_{b}$.

Since different D-block/T-block subsequences might have different landing/takeoff time increments, we can search for the optimal solutions for the optimization problem (\ref{optim1}) {by increasing or decreasing} the number, the landing/takeoff time increments and the time lengths of D-block/T-block subsequences of the sequence $\Phi_{c4}$.
Note that the quantity $F(\Phi_{c4})-F(\Phi_a)-F(\Phi_b)$ might correspond to several combinations of the landing/takeoff time increments and the increments of the time lengths of different groups of D-block/T-block subsequences compared to benchmark landing/takeoff sequences $\Phi_a$ and $\Phi_b$. Then we can obtain all the possible combinations of the landing/takeoff time increments and the increments of the time lengths of D-block/T-block subsequences such that the increments of the total delays are no larger than the quantity $F(\Phi_{c4})-F(\Phi_a)-F(\Phi_b)$. Then based on the obtained combinations, we can search for  the optimal solutions for the optimization problem (\ref{optim1}) within $[F(\Phi_a)+F(\Phi_b), \min\{F(\Phi_{c3}), F(\Phi_{d3})\}]$.

To implement the full space search in this step, we can fix the orders of the landing/takeoff aircraft and make discussions based on the number of breakpoints.

Moreover, in this step, we can use the parallel computing technology to reduce the computing time.

\begin{remark}{\rm The main idea of Algorithm 2 is to decompose the aircraft sequence into several block subsequences and fully explore combinations of the block subsequences along an optimal landing/takeoff subsequence to consider the optimization problem (\ref{optim1}). Note that the number of all possible block subsequences of different class combinations is not large according to the minimum separation time standards of Heathrow Airport and the RECAT system, when resident-point aircraft are not taken into account. The use of block subsequences might significantly reduce the computation amount of the algorithm.}\end{remark}

\begin{table}[!tbp]
		\caption{Minimum landing separation times in Heathrow Airport (Sec)} 
		\label{tab:Heathrow_separation}
		\centering
        \scriptsize
		\begin{tabular}{ccccccccccccc} 
			\toprule 
			\midrule 
            \multicolumn{2}{c}{} & \multicolumn{6}{c}{Trailing Aircraft} \\
            \midrule 
            \multicolumn{2}{c}{} & A & B & C & D & E & F \\
            \midrule 
            \multirow{6}*{Leading Aircraft}  & A & 90 & 135& 158& 158 & 158 & 180 \\
                                 & B & 90 & 90 & 113& 113 & 135 & 158 \\
                                 & C & 60 & 60 & 68 & 90  & 90  & 135 \\
                                 & D & 60 & 60 & 60 & 60  & 68  & 113 \\
                                 & E & 60 & 60 & 60 & 60  & 68  & 90 \\
                                 & F & 60 & 60 & 60 & 60  & 60  & 60 \\

			\midrule 
			\bottomrule 
		\end{tabular}
\end{table}
\begin{table}[!tbp]
		\caption{Minimum takeoff separation times based on RECAT-EU (Sec)} 
		\label{tab:Taking-off_separation}
		\centering
        \scriptsize
		\begin{tabular}{ccccccccccccc} 
			\toprule 
			\midrule 
            \multicolumn{2}{c}{} & \multicolumn{6}{c}{Trailing Aircraft} \\
            \midrule 
            \multicolumn{2}{c}{} & A & B & C & D & E & F \\
            \midrule 
            \multirow{6}*{Leading Aircraft}  & A & 80 & 100& 120& 140 & 160 & 180 \\
                                 & B & 80 & 80 & 100& 100 & 120 & 140 \\
                                 & C & 60 & 60 & 80 & 80  & 100 & 120 \\
                                 & D & 60 & 60 & 60 & 60  & 60  & 120 \\
                                 & E & 60 & 60 & 60 & 60  & 60  & 100 \\
                                 & F & 60 & 60 & 60 & 60  & 60  & 80 \\

			\midrule 
			\bottomrule 
		\end{tabular}
\end{table}

\section{Simulation}\label{simulations}
\begin{table*}[!htbp]
		\caption{{Comparison of performance and computation times for same-runway aircraft scheduling problem}} 
		\label{tab:cost_function CPU Time2}
		\centering
		\begin{tabular}{ccccccccccc} 
			\toprule 
			\midrule 
             \multirow{2}*{$T_W$(min)} & Aircraft Number  & \multirow{2}*{$T_E$(min)} & Operation & \multicolumn{2}{c}{Objective Function(s)} & \multicolumn{2}{c}{Computation Times(s)} & Gap\\
             & $|A|$  &  &  Task & MIP & Our Algorithm & MIP & Our Algorithm & \%\\
            \midrule 
            \multirow{12}*{60} & \multirow{3}*{30} & 20 & Takeoff & 14362 & 14362 & 600 & 0.20 & 0 \\
                              & & 20 & Landing & 8763 & 8279 & 600 & 0.15 & 5.85 \\
                              & & 20 & Mixed & 8671 & 8569 & 600 & 1.28 & 1.19\\
            \cmidrule{2-9}
             & \multirow{3}*{40} & 20 & Takeoff & 36695 & 35095 & 600 & 0.55 & 4.56  \\
                              &  & 20 & Landing & 28297 & 25234 & 600 & 0.37 & 12.14  \\
                              &  & 20 & Mixed   & 24520 & 24400 & 600 & 1.93 & 0.49\\
            \cmidrule{2-9}
             &\multirow{3}*{50} & 30 & Takeoff & 47607 & 45813 & 600 & 1.27 & 3.92 \\
                            &   & 30 & Landing & 37264 & 33141 & 600 & 1.11 & 12.44 \\
                            &   & 20 & Mixed   & 41914 & 41490 & 600 & 2.86 & 1.02 \\
            \cmidrule{2-9}
             &\multirow{3}*{60} & 30 & Takeoff & 77118 & 73478 & 600 & 3.24 & 4.95 \\
                           &    & 30 & Landing & 68173 & 58214 & 600 & 2.69 & 17.11 \\
                           &    & 20 & Mixed   & 75975 & 72852 & 600 & 4.56 & 4.29 \\
			\midrule 
            \multirow{12}*{90} & \multirow{3}*{30} & 20 & Takeoff & 12164 & 12144 & 600 & 0.18 & 0.16 \\
                              & & 20 & Landing & 9073 & 8785 & 600 & 0.16 & 3.29 \\
                              & & 20 & Mixed & 8453 & 8439 & 600 & 1.21 & 0.17\\
            \cmidrule{2-9}
             & \multirow{3}*{40}& 20 & Takeoff & 35835 & 35095 & 600 & 0.55 & 2.11  \\
                              & & 20 & Landing & 22738 & 22566 & 600 & 0.51 & 0.76  \\
                              & & 20 & Mixed & 24700 & 24400 & 600 & 1.95 & 1.23\\
            \cmidrule{2-9}
             &\multirow{3}*{50} & 20 & Takeoff & 61584 & 60164 & 600 & 1.35 & 2.36 \\
                            &   & 20 & Landing & 47067 & 42097 & 600 & 1.42 & 11.81 \\
                            &   & 20 & Mixed   & 43788 & 43424 & 600 & 2.87 & 0.84\\
            \cmidrule{2-9}
             &\multirow{3}*{60} & 20 & Takeoff & 89136 & 84736 & 600 & 2.65 & 5.19 \\
                           &    & 20 & Landing & 79896 & 71478 & 600 & 3.19 & 11.78 \\
                           &    & 20 & Mixed   & 67749 & 66484 & 600 & 4.32 & 1.90 \\
            \midrule 
            \multirow{12}*{120} & \multirow{3}*{30} & 20 & Takeoff & 15637 & 15637 & 600 & 0.15 & 0 \\
                            &   & 20 & Landing & 9120 & 8928 & 600 & 0.14 & 2.15\\
                            &   & 20 & Mixed & 6128 & 6109 & 600 & 1.21 & 0.31\\
            \cmidrule{2-9}
             &\multirow{3}*{40} & 20 & Takeoff & 28409 & 28069 & 600 & 0.57 & 1.21  \\
                            &   & 20 & Landing & 24201 & 21990 & 600 & 0.48 & 10.05  \\
                            &   & 20 & Mixed & 20906 & 20780 & 600 & 2.09 & 0.61\\
            \cmidrule{2-9}
             &\multirow{3}*{50} & 20 & Takeoff & 56504 & 53864 & 600 & 1.37 & 4.90  \\
                            &   & 20 & Landing & 53176  & 46422 & 600 & 1.29 & 14.55 \\
                           &    & 20 & Mixed & 37394 & 36829 & 600 & 3.19 & 1.53\\
            \cmidrule{2-9}
             &\multirow{3}*{60} & 20 & Takeoff & 86716 & 81845 & 600 & 2.89 & 5.95 \\
                         &      & 20 & Landing & 75769 & 68811 & 600 & 2.62 & 10.11 \\
                          &     & 20 & Mixed & 72386 & 71109 & 600 & 4.61 & 1.80\\
			\bottomrule 
		\end{tabular}
\end{table*}

\begin{table*}[htbp]
		\caption{{Comparison of performance and computation times for dual-runway aircraft scheduling problem}} 
		\label{tab:cost_function CPU Time3}
		\centering
		\begin{tabular}{ccccccccc} 
			\toprule 
			\midrule 
            $T_W$ & $T_E$ & Aircraft Number  & \multicolumn{2}{c}{Objective Function(s)} & \multicolumn{2}{c}{Computation Times(s)} & Gap\\
            (min)& (min) & $|A|$   & MIP & Our Algorithm & MIP & Our Algorithm & \%\\
            \midrule 
            \multirow{4}*{60} & 20& 70  & 47803 & 43066 & 600 & 4.15 & 11.00 \\
             &20& 80  & 71524 & 60823 & 600 & 4.57 & 17.59  \\
             &20& 90  & 88989 & 76128 & 600 & 6.77 & 16.89  \\
             &20& 100 & 131632 & 115304 & 600 & 8.79 & 14.16  \\
			\midrule 
            \multirow{4}*{90} & 20& 70  & 53389 & 46611 & 600 & 4.04 & 14.54 \\
             &30& 80  & 63372 & 59127 & 600 & 5.40 & 7.18  \\
             &40& 90  & 94272 & 82905 & 600 & 6.83 & 13.71  \\
             &50& 100 & 134749 & 112959 & 600 & 9.40 & 19.29  \\
            \midrule 
            \multirow{4}*{120} & 20& 70  & 47169 & 43171 & 600 & 3.66 & 9.26 \\
             &20& 80  & 74232 & 63237 & 600 & 5.57 & 17.39  \\
             &20& 90  & 96782 & 82216 & 600 & 8.16 & 17.71  \\
             &20& 100 & 121368 & 101823 & 600 & 9.57 & 19.20  \\
			\bottomrule 
		\end{tabular}
\end{table*}

In this section, we evaluate the efficiency of the proposed algorithm by comparing its computation time and objective function performance against a standard MIP solver. Two typical aircraft scheduling problems are considered: operations on a same runway and on dual runways with spacing no larger than 760 m. For the same-runway case, we examine four scenarios with $|A| = 30, 40, 50, 60$ under three different operations: takeoff only, landing only, and mixed operations (both takeoff and landing). For the dual-runway case, we consider scenarios with $|A| = 70, 80, 90, 100$ aircraft, each comprising an equal number of takeoff aircraft and landing aircraft. The aircraft are classified into six categories based on the RECAT-EU framework $(A, B, C, D, E, F)$, with proportions of $10\%$, $20\%$, $25\%$, $15\%$, $20\%$, and $10\%$, respectively.

The minimum separation times for aircraft pairs with the same type of operation are provided in Tables~\ref{tab:Heathrow_separation} and \ref{tab:Taking-off_separation}. For both problems, the separation time between a takeoff aircraft and a following landing aircraft is set to $D_T = D_P = 60$ seconds, while the separation time between a landing aircraft and a following takeoff aircraft is set to $T_D = 75$ seconds for the same-runway case and $P_D = 0$ seconds for the dual-runway case.

All simulations are conducted on a computer equipped with an AMD Ryzen 7 7840H processor (3.8 GHz, 16 GB RAM). The MIP formulations are solved using CPLEX 12.10 with a time limit of $600$ seconds per instance. The proposed algorithm is implemented in MATLAB R2021b. To enhance computational efficiency, we employ MATLAB's parallel computing technology using $8$ workers to deal with the mixed operation scenarios on a same runway, and the dual-runway scheduling scenarios.

For the same-runway case, the earliest landing or takeoff times for each aircraft are randomly generated, following a uniform distribution within the interval $[0, T_E]$ minutes. The time window lengths are set to $T_W = 60$, $90$, and $120$ minutes to evaluate the algorithms' performance under different levels of scheduling flexibility. The scheduled time for each aircraft is set as $P_i = f_i^{\text{min}}$ for takeoff aircraft and $P_i = f_i^{\text{min}} + 5$ minutes for landing aircraft. Table~\ref{tab:cost_function CPU Time2} presents a comparison of the objective function values and computation times between the proposed algorithm and the MIP solver, along with the percentage gap in objective values between the two methods.

For the dual-runway case, aircraft earliest operation times are also uniformly distributed over $[0, T_E]$ minutes. The time window lengths remain at $T_W = 60$, $90$, and $120$ minutes.
Table~\ref{tab:cost_function CPU Time3} summarizes the performance comparison between our algorithm and the MIP solver in terms of objective function values and computation time, as well as the corresponding objective value gaps.

From Tables~\ref{tab:cost_function CPU Time2} and \ref{tab:cost_function CPU Time3}, it is evident that the proposed algorithm consistently obtains solutions within $10$ seconds, whereas the MIP solver requires the full $10$-minute time limit. Moreover, our algorithm consistently yields better objective function values than those produced by the MIP solver, with the performance gap widening as the number of aircraft increases. In addition, the runtime of the MIP solver has been extended to over one hour but no significant improvements were found in the objective function values compared to the results obtained within the $10$-minute time limit.

In the same-runway scenario, our algorithm shows particularly strong performance in the takeoff-only and landing-only cases. The performance advantage is less pronounced in the mixed takeoff-and-landing case. This is primarily due to the relatively short separation times required between different operation tasks: $60$ seconds between a takeoff aircraft and a trailing landing aircraft, and $75$ seconds between a landing aircraft and a trailing takeoff aircraft. These moderate separation values enable the MIP solver to generate relative high-quality solutions by frequently using takeoff-landing and landing-takeoff transitions. Nonetheless, as the number of aircraft increases, the superiority of our algorithm remains increasingly evident.

For the dual-runway scenario, we evaluate $12$ instances with progressively larger numbers of aircraft, reflecting the higher operational capacity of dual-runway airports. In particular, when the problem size reaches $90$ and $100$ aircraft, the performance gap in objective function values between our algorithm and the MIP solver exceeds $17\%$, exhibiting the scalability and efficiency of our proposed algorithms.

\section{Conclusions}

In this paper, scheduling problems of aircraft minimizing the total delays on a same runway and on dual runways are studied.
Two real-time optimal algorithms were proposed for the four scheduling problems by fully exploiting the combinations of different classes of aircraft based on parallel computing technology.
When $100$ aircraft on dual runways were considered, by using the algorithm in this paper, the optimal solution can be obtained within less than $10$ seconds, while by using the CPLEX software to solve the mix-integer optimization model, the optimal solution cannot be obtained within $1$ hour.

\end{document}